\documentclass[11pt,reqno]{amsart}
                                                                                
\usepackage{graphicx}
\usepackage{amscd}
\usepackage{dsfont}
\usepackage{amsmath}
\usepackage{amsfonts}
\usepackage{psfrag}
\usepackage{amssymb}
\usepackage{verbatim}
\usepackage{comment}
\newtheorem{theorem}{Theorem}[section]
\theoremstyle{plain}

\newtheorem{cor}[theorem]{Corollary}

\newtheorem{definition}[theorem]{Definition}

\newtheorem{lemma}[theorem]{Lemma}

\newtheorem{prop}[theorem]{Proposition}
\newtheorem{remark}[theorem]{Remark}

\numberwithin{equation}{section}

\newcommand{\lam}{\lambda}
\newcommand{\gam}{\gamma}
\newcommand{\om}{\omega}
\newcommand{\Sig}{\Sigma}

\newcommand{\Gam}{\Gamma}
\newcommand{\sig}{\sigma}
\newcommand{\by}{{\bf y}}
\newcommand{\bi}{{\bf i}}
\newcommand{\bj}{{\bf j}}
\newcommand{\R}{{\mathbb R}}

\newcommand{\E}{{\mathbb E}\,}
\newcommand{\Prob}{{\mathbb P}\,}
\newcommand{\Nat}{{\mathbb N}}
\newcommand{\Leb}{{\mathcal L}}
\newcommand{\Ek}{{\mathcal E}}

\newcommand{\Wk}{{\mathcal W}}
\newcommand{\eps}{{\varepsilon}}
\newcommand{\es}{\emptyset}
\newcommand{\udim}{\overline{\dim}}
\newcommand{\const}{{\rm const}}
\newcommand{\Dh}{{\dim_{_{\scriptstyle H}}}}
\newcommand{\nuhat}{\widehat{\nu}}

\def\ess{{\rm ess}}
\def\lt{\left}                \def\rt{\right}
\def\Sk{{\mathcal S}}
\def\half{\frac{1}{2}}

\begin{document}
\title[IFS with random maps]
{Absolute continuity for random iterated function systems with overlaps}

\author{Yuval peres}
\address{Yuval Peres, Department of Statistics, University of
California, Berkeley} \email{peres@stat.Berkeley.edu}

\author{K\'{a}roly Simon}
\address{K\'{a}roly Simon, Institute of Mathematics, Technical
University of Budapest, H-1529 B.O.box 91, Hungary}
\email{simonk@math.bme.hu}

\author{Boris Solomyak}
\address{Boris Solomyak, Box 354350, Department of Mathematics,
University of Washington, Seattle WA 98195}
\email{solomyak@math.washington.edu}

 \thanks{
\indent Research of Peres was
partially supported by NSF grants \#DMS-0104073 and \#DMS-0244479.
Part of this work was done while he was visiting Microsoft Research.
Research of Solomyak was partially supported by NSF grants
 \#DMS-0099814 and \#DMS-0355187.
Research of
Simon was partially supported by OTKA Foundation grant \#T42496.
The collaboration of K.S. and B.S. was supported by NSF-MTA-OTKA grant \#77.} 

\begin{abstract}
We consider linear iterated function systems with a random multiplicative
error on the real line. Our system is $\{x\mapsto d_i + \lam_i Y x\}_{i=1}^m$,
where $d_i\in \R$ and $\lam_i>0$ are fixed and $Y> 0$
is a random variable with an absolutely continuous distribution.
The iterated maps are applied randomly
according to a stationary ergodic process,
with the sequence of i.i.d.\ errors $y_1,y_2,\ldots$, distributed as $Y$,
independent of everything else. Let $h$ be the entropy of the process,
and let $\chi = \E[\log(\lam Y)]$ be the Lyapunov
exponent. Assuming that $\chi < 0$,
we obtain a family of conditional measures $\nu_\by$ on the line,
parametrized by $\by = (y_1,y_2,\ldots)$, the sequence of errors.
Our main result is that if $h > |\chi|$, then $\nu_\by$ is absolutely
continuous with respect to the Lebesgue measure for a.e.\ $\by$.
We also prove that if $h < |\chi|$, then
the measure $\nu_\by$ is singular and has dimension $h/|\chi|$ for
a.e.\ $\by$. These results are applied to a randomly perturbed IFS 
suggested by Y. Sinai, and
to a class of random sets considered by R. Arratia, motivated by
probabilistic number theory.
\end{abstract}
\date{\today}

\maketitle

\thispagestyle{empty}

\section{Introduction} 

Let $\{f_1,\ldots,f_m\}$ be an iterated function system (IFS) on the real line,
where the maps are applied according to the probabilities $(p_1,\ldots,p_m)$,
with the choice of the map random and independent at each step.
We assume that the system is contracting on average,
that is, the Lyapunov exponent $\chi$ (appropriately defined) is
negative. In this paper the maps will be linear, $f_i(x)=\lam_i x+ d_i$, and
then $\chi:= \sum_{i=1}^m p_i \log \lam_i$. If $\chi<0$, then there is a
well-defined invariant probability measure $\nu$ on $\R$ (see \cite{DiFr}).
It is of interest to determine whether this measure is singular or absolutely
continuous, and if it is singular, to compute its Hausdorff dimension
$$
\Dh (\nu) = \inf\{\Dh(Y):\ \nu(\R\setminus Y)=0\}.
$$
Let 
$h = -\sum_{i=1}^m p_i \log  p_i$ be the entropy of the underlying Bernoulli
process. It was proved in \cite{NSB} for non-linear 
contracting on average IFS (and later extended in \cite{FST}) that 
$$
\Dh (\nu) \le h/|\chi|.
$$
A question arises what happens when the entropy is greater than the 
absolute value of the Lyapunov exponent. One can expect that, at least 
``typically,'' the measure $\nu$ is absolutely continuous when $h/|\chi|>1$.
Such results are known for contracting IFS (see \cite{peso2,SSU2,Neu,ngawang}),
but extending them
to the case when it is only contracting on average remains a challenge.

In this paper we study a modification of the problem which makes it more
tractable, namely we consider linear IFS with a random multiplicative error.
Our system is $\{x\mapsto d_i + \lam_i Y x\}_{i=1}^m$,
where $d_i\in \R$ and $\lam_i>0$ are fixed and $Y\ge 0$
is a random variable with an absolutely continuous distribution.
The iterated maps are applied randomly
according to a stationary ergodic process,
with the sequence of i.i.d.\ errors $y_1,y_2,\ldots$, distributed as $Y$,
independent of everything else. Let $h$ be the entropy of the process, 
and let $\chi = \E[\log(\lam Y)]$ be the Lyapunov
exponent (the symbol $\E$ denotes expectation). Assuming that $\chi < 0$, 
we obtain a family of conditional measures $\nu_\by$ on the line,
parametrized by $\by = (y_1,y_2,\ldots)$, the sequence of errors. 
Our main result is that if $h > |\chi|$, then $\nu_\by$ is absolutely
continuous with respect to the Lebesgue measure $\Leb$ for a.e.\ $\by$.
We also prove that if $h < |\chi|$, then
the measure $\nu_\by$ is singular and has dimension $h/|\chi|$ for
a.e.\ $\by$. Random IFS are quite well understood under separation 
conditions (see \cite{Ki1,Ki2}), and the novelty here is that overlaps are
allowed.

Before stating the general results precisely,
we describe two examples that motivated our work. Consider the random series
\begin{equation} \label{sumprod}
X = 1 + Z_1 + Z_1 Z_2 + \ldots + Z_1 Z_2\cdots Z_n + \ldots
\end{equation}
where $Z_i$ are i.i.d.\ random variables, taking the values in 
$\{1-a,1+a\}$ with probabilities $(\half,\half)$, for a fixed parameter 
$a\in (0,1)$. 
The series converges almost surely, since $\chi=
\E[\log Z] = \half\log(1-a^2)<0$.
Let $\nu^a$ denote the distribution of $X$. 
Sinai [personal communication], motivated by
a statistical analog of the well-known open ``$3n+1$ problem,''
asked for which $a$ the distribution of $X$ is absolutely continuous.
Observe that $\nu^a$ is supported on $[\frac{1}{a},\infty)$ for all $a>0$.
Note that $\nu^a$ is the invariant measure for the IFS $\{1+(1-a)x,
1+(1+a)x\}$, with probabilities $(\half,\half)$. 
If $h < |\chi|$, that is, $\log 2 < -\half \log(1-a^2)$, or
$a > \sqrt{3}/2$, then the measure is singular.
It is natural to predict that $\nu^a \ll \Leb$ for a.e.\ $a\in (0,\sqrt{3}/2)$.
While this question remains open, here we solve a ``randomly perturbed
version.'' 

\begin{prop} \label{prop-ex1} Consider the random sum (\ref{sumprod}), 
with $Z_i = \lam_i Y$, where $\lam_i \in \{1-a,1+a\}$ with probabilities 
$(\half,\half)$ and $Y$ has an absolutely continuous distribution on 
$(1-\eps_1,1+\eps_2)$ for small $\eps_1$ and $\eps_2$,
with a bounded density, such that $\E[\log Y] = 0$.  The ``errors''
$y_i$ at each step are i.i.d.\ with the distribution of $Y$, and are
independent of everything else. Let $\nu^a_\by$ be the conditional distribution
of $Z$, given a sequence of errors $\by = (y_1,y_2,\ldots)$.

{\bf (a)} If $a\in (0,\sqrt{3}/2)$, then   $\nu^a_\by \ll \Leb$
for a.e.\ $\by$;

{\bf (b)} if $a\ge \sqrt{3}/2$, then $\nu^a_\by \perp \Leb$ and
$\Dh (\nu^a_\by) = 2\log 2/\log\frac{1}{1-a^2}$ for a.e.\ $\by$.
\end{prop}

Our second example is probabilistic in its origin (rather than a random
perturbation of a deterministic one, as above). It comes from a question of
Arratia (see Section 22 in \cite{arratia}),
who considered the following distributions, motivated by
some questions in probabilistic number theory.
Let $Y = U^{1/\theta}$, where $U$ has the uniform distribution on $[0,1]$,
and consider $X_i = Y_1\cdots Y_i$, where $Y_i$ are i.i.d.\ with the
distribution of $Y$. The process $\{X_i\}$ is known as 
the scale-invariant Poisson process with intensity $\theta x^{-1}\, dx$.
Consider the random sum $Z = \sum_{i\ge 1}
J_i X_i$ where $J_1, J_2,\ldots $ are ``fair coins'' with values in $\{0,1\}$,
independent of each other and of everything else.
One is interested in the conditional distribution $\nu^\theta_\by$
of $Z$ given the process $\{Y_i\}$, and in its support,
$S_\by^\theta = \{\sum_{i=1}^\infty a_i X_i:\ a_i \in \{0,1\}\}$.
Observe that this fits into our class of IFS $\{d_i+\lam_i Yx\}_{i=1}^m$,
by taking $m=2$, $d_1=0, d_2=1, \lam_1=\lam_2=1$.
The distribution of $ U^{1/\theta}$ has
the density $\theta x^{\theta-1} \mathds{1}_{[0,1]}$, so we have
$\chi=\E[\log \lam Y] = \E[\log Y] = -\theta^{-1}$.
The entropy of the ``fair coins''
process is  $h=\log 2$, so $h/|\chi| = \theta \log 2$.

\begin{prop} \label{prop-ex2}
Let $Z = \sum_{i\ge 1}
J_i X_i$ as above.
Let $\nu^\theta_\by$ be the conditional distribution of $Z$ given the 
process $\{Y_i\}$, and let $S_\by^\theta$ be 
the support of $\nu^\theta_\by$.

{\bf (a)} If $\theta > 1/\log 2$, then $\nu^\theta_\by\ll \Leb$ with a 
density in $L^2(\R)$, hence $\Leb(S_\by^\theta)>0$, for a.e.\ $\by$.
                                                                                
{\bf (b)} If $\theta > 2/\log 2$, then $\nu^\theta_\by\ll \Leb$ with a
continuous density, hence $S_\by^\theta$ contains an interval,
for a.e.\ $\by$.

{\bf (c)} If $\theta \in (0,1/\log 2]$, then
$\Leb(S_\by^\theta)=0$, 
and $\Dh(S_\by^\theta)=\Dh (\nu^\theta_\by) = \theta \log 2$
for a.e.\ $\by$.
\end{prop}

An intriguing open problem is whether $S_\by^\theta$ contains intervals
for a.e.\ $\by$ when $\theta \in (\frac{1}{\log 2},\frac{2}{\log 2})$.
The proof of Proposition~\ref{prop-ex2} 
is easily adapted to the Poisson-Dirichlet
distributions where the variables $X_i = Y_1\cdots Y_i$ are replaced by
$\widetilde{X}_i = Y_1\cdots Y_{i-1} (1-Y_i)$ ordered by size, see
the equations (6.2) and (8.2) in \cite{arratia}.

\begin{sloppypar}
Our results yield many variants of Proposition~\ref{prop-ex2}.
For example, let $\widetilde{S}_\by^\theta = 
\{\sum_{i=1}^\infty a_i X_i:\ a_i \in \{0,1\},\,a_i a_{i+1}=0, i\ge 1\}$;
in other words, we consider only the sums corresponding to the ``Fibonacci''
shift of finite type. Let $\tau = (1+\sqrt{5})/2$.
\end{sloppypar}

\begin{prop} \label{prop-ex2a}
{\bf (a)} If $\theta > 1/\log \tau$, then 
$\Leb(\widetilde{S}_\by^\theta)>0$, for a.e.\ $\by$.

{\bf (b)} If $\theta \in (0,1/\log \tau)$, then
$\Leb(\widetilde{S}_\by^\theta)=0$ and $\Dh(\widetilde{S}_\by^\theta)=
\theta \log\tau$ for a.e.\ $\by$.
\end{prop}
 

\section{Statement of results}

Consider a random variable $Y$ with an absolutely continuous
distribution $\eta$ on $(0,\infty)$, such that for some $C_1>0$ we have
\begin{equation} \label{eq-density}
\frac{d\eta}{d x} \le C_1 x^{-1},\ 
\forall\,x > 0.
\end{equation}
Let $\R^{\Nat}$ be the
infinite product equipped with the product measure 
$\eta_{\infty }:=\eta^\Nat$. Let $\mu$ be an
ergodic $\sigma$-invariant measure on 
$\Sigma=\left\{1,\dots ,m\right\}^{\mathbb{N}}$, where
$\sigma$ is the left shift. Denote by $h(\mu)$ the entropy of the measure
$\mu$. We consider linear IFS with a random multiplicative error
$\{x\mapsto d_i + \lam_i Y x\}_{i=1}^m$,
where $d_i\in \R$ and $\lam_i>0$ are fixed.
The iterated maps are applied randomly
according to the stationary measure $\mu$,
with the sequence of i.i.d.\ errors $y_1,y_2,\ldots$, distributed as $Y$,
independent of the choice of the function.
The Lyapunov exponent of the IFS is defined by
$$
\chi(\mu,\eta):= \E[\log (\lam Y)] = \E[\log Y] + 
\int_\Sigma \log(\lam_{i_1})\,d\mu(\bi).
$$
Throughout the paper, we assume that
\begin{equation} \label{contract}
\chi(\mu,\eta) < 0,
\end{equation}
which means that the IFS is contracting on average.
The natural projection $\Pi :\Sigma
\times \R^\Nat\to \R$ is defined by
\begin{equation} \label{eq-proj}
\Pi (\mathbf{i},\mathbf{y}):=%
d_{i_1}+d_{i_2}\lambda_{i_1}y_1+\cdots +
d_{i_{n+1}}\lambda_{i_1\dots i_n}y_{1 \ldots n}+\ldots,
\end{equation}
where $\bi = (i_1,i_2,\ldots)$,  
$\lambda_{i_1\dots i_n}:=\lambda_{i_1}\cdots \lambda_{i_n}$, and
$y_{1 \ldots n}=y_1\cdots y_n$. 
Note that $\Pi(\bi,\by)$ is a Borel map
defined $\mu\times \eta_\infty$ a.e., since 
$
n^{-1} \log (\lam_{1\ldots n} y_{1\ldots n}) \to \E[\log (\lam Y)]<0
$
a.e., by the Birkhoff Ergodic Theorem.
For a fixed $\by\in \R^\Nat$ we define 
$\Pi_{\by}:\Sigma \to \R$ and the measure $\nu_{\by}$ on $\R$ as
\begin{equation} \label{eq-nuby}
\Pi_{\by}(\bi ):=\Pi (\bi,\by) \ \ \mbox{ and } \ \ 
\nu_{\by}:=\left(\Pi_{\by}\right)_{*}\mu.
\end{equation}

We need to impose a condition which guarantees that the maps of the IFS
are sufficiently 
different. We consider two cases which cover the interesting examples
that we know of. We assume that either all the digits are distinct:
\begin{equation} \label{cond1}
d_i \ne d_j,\ \ \ \mbox{for all}\ i\ne j,
\end{equation}
or all the digits $d_i$ are equal to some $d\ne 0$ (which we can assume to
be 1, without loss of generality),
but the average contraction ratios are all distinct: 
\begin{equation} \label{cond2}
d_i = 1,\ \ \lam_i \ne \lam_j,\ \ \ \mbox{for all}\ i\ne j.
\end{equation}

\begin{theorem} \label{main}
Let $\nu_\by$ be the conditional distribution of the sum (\ref{eq-proj})
given $\by=(y_1,y_2,\ldots)$. We assume that
(\ref{eq-density}), (\ref{contract}) hold, and either (\ref{cond1}) or
(\ref{cond2}) is satisfied.

{\bf (a)} 
If $h(\mu) > |\chi(\mu,\eta)|$, then
\begin{equation} \label{eq-main}
\nu_{\by} \ll \Leb \ \ \mbox{for}\ \eta_\infty\ \mbox{a.e.\ } \by.
\end{equation}

{\bf (b)} 
If $h(\mu) \le |\chi(\mu,\eta)|$, then 
\begin{equation} \label{meas-dim}
\Dh (\nu_\by) = \frac{h(\mu)}{|\chi(\mu,\eta)|} \ \ \mbox{for}\ 
\eta_\infty\ \mbox{a.e.\ } \by.
\end{equation}
\end{theorem}

Assuming that $\mu$ is a product (Bernoulli) measure, that is,
$\mu = (p_1,\ldots,p_m)^\Nat$ and $h(\mu)=|\chi(\mu,\eta)|$, we can show 
that the measure $\nu_\by$ is singular for a.e.\ $\by$.

\begin{prop} \label{prop-sing} Suppose that $\mu = (p_1,\ldots,p_m)^\Nat$
and $Y>0$ is any random variable, such that (\ref{contract}) holds and
and $(\bi,\by) \mapsto p_{i_1} \lam_{i_1}^{-1} y_1^{-1}$ is 
non-constant on $\Sig\times \R^\Nat$. If $h(\mu)=-\chi(\mu,\eta)$,
then $\nu_\by \perp \Leb$ for $\eta_\infty$ a.e.\ $\by$.
\end{prop}

Note that in the proposition we do not make any assumptions on the
distribution of $Y$. If $Y$ is any non-constant
random variable, then the proposition applies.
On the other hand, it includes the case when $Y$ is
constant (in other words, this is a usual IFS with no randomness), but
$\lam_j/p_j$ is not constant. Then of course $Y$ can be eliminated
altogether and there is no a.e.\ $\by$ in the statement.
In the special case $Y\equiv 1$ and $\lam_i<1$ for all $i\le m$, our statement 
is contained in \cite[Th.1.1(ii)]{ngawang}. 

We should emphasize that Proposition~\ref{prop-sing}, as well as
the upper dimension estimate in (\ref{meas-dim}), are rather
standard; they are
included for completeness, in order to indicate that our results are sharp.

Next we discuss two special cases which include the examples from the
Introduction.

\subsection{Sums of products of i.i.d.\ random variables}
Suppose that (\ref{cond2}) holds and 
$\mu=(p_1,\ldots,p_m)^\Nat$. 
Then we are led to a random variable
\begin{equation} \label{sumpro}
X = 1 + Z_1 + Z_1 Z_2 + \ldots + Z_1 Z_2\cdots Z_n + \ldots
\end{equation}
where $Z_i$ are independent 
with the distribution of 
$Z=\lam Y$, where $\lam = \lam_j$ with probability $p_j$
and $Y$, independent of $\lam$, satisfies (\ref{eq-density}), as
in the general case.
Then $\chi(\mu,\eta) = \sum_{i=1}^m p_i \log \lam_i+ \E[\log Y]$ and
$h(\mu) = -\sum_{i=1}^m p_i \log p_i$.
The measure $\nu_\by$ is the conditional
distribution of $X$ given $\by = (y_1,y_2,\ldots)$, a realization of the
process $\{Y_i\}$.

Thus Proposition~\ref{prop-ex1} on the randomly perturbed Sinai's problem
is a special case of Theorem~\ref{main} and Proposition~\ref{prop-sing}.

\subsection{Homogeneous case: random measures} \label{subs-Arr}

Suppose that $\lam_i=\lam$ for all $i\le m$, so we have $d_i\ne d_j$ for
$i\ne j$ by (\ref{cond1}). 
We are led to the random sums 
\begin{equation} \label{absorb}
X = \sum_{i=1}^\infty J_{i} Y_1\cdots Y_{i-1}\lam^{i-1},
\end{equation}
where $Y_i$ are i.i.d.\  with the absolutely continuous distribution $\eta$,
and $J_i$ take the values in $\{d_1,\ldots,d_m\}$,
are independent of $\{Y_j\}$, and are chosen
according to an ergodic $\sig$-invariant measure $\mu$.
Then Theorem~\ref{main} applies to $\nu_\by$, the conditional
distribution of $X$ given $\by = (y_1,y_2,\ldots)$, a realization of the
process $\{Y_j\}$. The Lyapunov exponent $\chi(\eta) = \log\lam+\E[\log Y]$ 
does not depend on $\mu$. 

When $\mu$ is Bernoulli, that is, $\mu = (p_1,\ldots,p_m)^\Nat$, there is
an alternative method to study $\nu_\by$ which goes back to the work of
Kahane and Salem \cite{KS} and uses Fourier transform. It requires a
stronger assumption on the distribution $\eta$, namely that
\begin{equation} \label{eq-dens2}
\eta\ \ \mbox{has compact support and}\ \ d\eta/dx\ \ \mbox{is of bounded
variation}.
\end{equation}

\begin{theorem} \label{th-homog}
Let $\nu_\by$ be the conditional distribution of the sum (\ref{eq-proj})
given $\by=(y_1,y_2,\ldots)$, defined by (\ref{eq-nuby}). We assume that
(\ref{eq-dens2}) and (\ref{cond1}) are
satisfied, $\mu = (p_1,\ldots,p_m)^\Nat$
and $\lam_i=\lam$ for all $i\le m$.
Suppose that $\chi(\eta) = \log\lam+\E[\log Y]<0$.

{\bf (a)} If $|\log(\sum_{i=1}^m p_i^2)| > |\chi(\eta)|$, 
then $\nu_\by \ll \Leb$ with a density in $L^2(\R)$ for $\eta_\infty$
a.e.\ $\by$.
                                                                                
{\bf (b)} If $|\log(\sum_{i=1}^m p_i^2)| > 2 |\chi(\eta)|$, 
then $\nu_\by \ll \Leb$ with a continuous density for $\eta_\infty$ a.e.\ $\by$.
\end{theorem}

Observe that in the uniform case, when $p_i = \frac{1}{m}$ 
for $i\le m$, we get the
same threshold $|\chi(\eta)|=\log m$ 
for absolute continuity in Theorem~\ref{main}(a)
and for absolute continuity with a density in $L^2$,
in Theorem~\ref{th-homog}(a).

\subsection{Homogeneous case: random sets}

The results on random measures yield information on random sets. 
Recall that $\Sig = \{1,\ldots,m\}^\Nat$, and
let $\Gam \subset \Sigma$ be a closed $\sigma$-invariant subset.
For a digit set $\{d_1,\ldots,d_m\}$ and $\by \in (0,\infty)^\Nat$ consider
$$
S_\Gam(\by) = \left\{\sum_{i=1}^\infty d_{a_i} y_{1\ldots (i-1)}:\
\{a_i\}_1^\infty \in \Gam \right\}\,.
$$ 
We let $S(\by) = S_\Sig(\by)$.
Denote by $h_{\rm top}(\Gam)$ the topological entropy of $(\Gam,\sig)$.
In the next three corollaries we consider a digit set
$\{d_1,\ldots,d_m\}$ satisfying (\ref{cond1}), that is, all the digits are
assumed to be distinct. For a random variable $\eta$ we let
$\chi(\eta):= \int \log \eta\,d\eta < 0$.

\begin{cor} \label{cor-ranset}
Suppose that $\eta$ satisfies (\ref{eq-density})
and $\chi(\eta) < 0$. 

{\bf (a)} If $h_{\rm top}(\Gam) > |\chi(\eta)|$, then $\Leb(S_\Gam(\by)) > 0$
for $\eta_\infty$ a.e.\ $\by$.

{\bf (b)} 
If $h_{\rm top}(\Gam) \le |\chi(\eta)|$, then 
$\Dh (S_\Gam(\by)) = h_{\rm top}(\Gam)/|\chi(\eta)|$ for $\eta_\infty$
a.e.\ $\by$.
\end{cor}

\begin{cor} \label{cor-ranset1}
Suppose that $\eta$ satisfies (\ref{eq-dens2})
and $\chi(\eta) < 0$.
We consider $\Gam = \Sig$.

{\bf (a)}  If $\log m > 2 |\chi(\eta)|$, then
$S(\by)$ contains an interval for $\eta_\infty$ a.e.\ $\by$.

{\bf (b)} If $\log m\le |\chi(\eta)|$, then 
$\Dh (S(\by)) = \log m/|\chi(\eta)|$ for $\eta_\infty$ a.e.\ $\by$.
\end{cor}

\begin{cor} \label{cor-zero}
Suppose that $Y$ is any non-constant random variable on $(0,\infty)$
such that $\chi(\eta) < 0$.
If $\log m = |\chi(\eta)|$, then
$\Leb(S(\by))=0$ for $\eta_\infty$ a.e.\ $\by$.
\end{cor}

The results of this subsection imply the statements on Arratia's question and
its variants from the Introduction. More precisely, parts (a) and (b) of
Proposition~\ref{prop-ex2} follow from Theorem~\ref{th-homog}. 
Part (c) of
Proposition~\ref{prop-ex2} follows from Theorem~\ref{main}(b) and 
Corollary~\ref{cor-zero}. Proposition~\ref{prop-ex2a} follows from
Corollary~\ref{cor-ranset}, since $\Gam = \{(a_i)_1^\infty \in 
\{0,1\}^\Nat:\ a_i a_{i+1} = 0, \,i\ge 1\}$  has topological entropy
$\log\frac{1+\sqrt{5}}{2}$.

The rest of the paper is organized as follows. Theorem~\ref{main}(a) is
proved in Sections 3 and 4; the latter also contains 
a key ``transversality lemma,''
which is used in the proof of both Theorem~\ref{main}(a) 
and the lower estimate in Theorem~\ref{main}(b).
Then Theorem~\ref{main}(b) is derived in Section 5.
Section 6 is devoted to the proofs of other results, especially
Theorem~\ref{th-homog}, following the method of Kahane and Salem.
Finally, Section 7 contains some open questions.


\section{Preliminaries and the Proof of Theorem~\ref{main}(a)}

{\em Notation.} For $\om\in \{1,\ldots,m\}^n$ 
we denote by $[\om]$ the cylinder set 
of $\bi \in \Sigma$
which start with $\om$. For $\bi \in \Sig$ let $[\bi,n] = [i_1\ldots i_n]$.
For $\bi,\bj\in \Sig$ we denote by $\bi\wedge \bj$ their common initial
segment.

By adding the constant $\E[\log Y]$ to $\log \lam$ and subtracting it from
$\log Y$, we can assume without loss of generality
that $\E[\log Y] = 0$, so that $\chi(\mu,\eta) = \chi(\mu)=
\E[\log \lam]$.  In order to prove Theorem~\ref{main},
we need to make a certain ``truncation'' both in $\R^\Nat$ and in $\Sigma$.
By the Law of Large Numbers,
\begin{equation} \label{eq-LLN}
n^{-1} \log(y_{1\ldots n}) \to 0\ \ \ \mbox{for}\ \eta_\infty\ \mbox{a.e.}
\ \by.
\end{equation}
By Egorov's Theorem, for any $\eps>0$ there exists
$F_\eps\subset \R^\Nat$, with $\eta_\infty(F_\eps)> 1-\eps$, such that
$(y_{1\ldots n})^{1/n} \to 1$ uniformly on $F_\eps$.

Next we do the truncation in $\Sigma$.
By the Shannon-McMillan-Breiman Theorem,
\begin{equation}\label{shannon}
n^{-1} \log(\mu[\bi,n]) \to -h(\mu)\ \ \ \mbox{for}\ \mu\ \mbox{a.e.}
\ \bi \in \Sigma.
\end{equation}
By the Birkhoff Ergodic Theorem,
\begin{equation}\label{ergth}
n^{-1} \log(\lam_{i_1\ldots i_n})\to \chi(\mu)
\ \ \ \mbox{for}\ \mu\ \mbox{a.e.} \ \bi \in \Sigma.
\end{equation}
Applying Egorov's Theorem,
we can find $G_\eps\subset \Sigma$, with
$\mu(G_\eps) > 1-\eps$, such that the convergence in (\ref{shannon}) and
(\ref{ergth}) is uniform on $G_\eps$.

Define $\mu_\eps = \mu |_{G_\eps}$ and let $\nu_\by^\eps = (\Pi_\by)_*\mu_\eps$.
We can work with measures $\nu_\by^\eps$ instead of $\nu_\by$. Indeed,
if $\nu_\by^\eps \ll \Leb$ for all $\eps>0$, then $\nu_\by \ll \Leb$, and
$\Dh (\nu_\by) = \sup_{\eps>0} \Dh (\nu_\by^\eps)$.
Since we can
obviously assume that $F_\eps \subset F_{\eps'}$ for $\eps'< \eps$, 
(\ref{eq-main}) will follow if we prove that
\begin{equation} \label{eq-main1}
\forall\,\eps>0,\ \nu_\by^\eps \ll \Leb\ \ \ \mbox{for}\ \eta_\infty
\ \mbox{a.e.}\
\by \in F_\eps.
\end{equation}
Similarly, (\ref{meas-dim}) will follow if we prove that
\begin{equation} \label{eq-main3}
\forall\,\eps>0,\ 
\Dh (\nu_\by^\eps) = \frac{h(\mu)}{|\chi(\mu)|} \ \ \mbox{for}\
\eta_\infty\ \mbox{a.e.\ } \by \in F_\eps.
\end{equation}

\medskip

\noindent
{\em Beginning of the Proof of Theorem~\ref{main}(a).}
Fix $\eps\in (0,1)$; our goal is to prove (\ref{eq-main1}) assuming that
$-h(\mu) < \chi(\mu) < 0$. We can fix positive $\theta < \rho$ such that
\begin{equation} \label{trunc3}
-h(\mu) < \log \theta  < \log \rho < \chi(\mu) < 0.
\end{equation}
Next fix $\delta > 0$ such that
\begin{equation}\label{delta}
(1+\delta) \theta < \rho.
\end{equation}
Using the uniform convergence on $F_\eps$ and $G_\eps$, we can find
$N = N(\eps,\delta,\theta,\rho)\in \Nat$ such that, in view of 
(\ref{eq-LLN}),
\begin{equation} \label{trunc1}
(1+\delta)^{-n} \le y_1\cdots y_n \le (1+\delta)^n\ \ \ \mbox{for all}\ \
n\ge N,\ \by \in F_\eps,
\end{equation}
and, in view of (\ref{trunc3}),
\begin{equation} \label{trunc4}
\mu[\bi,n]  < \theta^n < \rho^n  < \lam_{i_1\ldots i_n},
\ \ \ \mbox{for all}\ \
n\ge N,\ \bi \in G_\eps.
\end{equation}
We can decompose the measure into the  sum of measures on cylinders: 
\begin{equation} \label{decomp}
\nu^\eps_{\by} = \sum_{|\om| = N} \nu^\eps_{\by,\om},\ \ \ \mbox{where}\ \
\nu^\eps_{\by,\om}:= (\mu|_{[\om]\cap G_\eps})\circ \Pi_\by^{-1}.
\end{equation}
Thus it is enough to show that
\begin{equation} \label{eq-main2}
\nu^\eps_{\by,\om} \ll \Leb\ \ \ \mbox{for}\ \eta_\infty\ \mbox{a.e.}\
\by \in F_\eps,\ \forall\,\om,\ |\om|=N.
\end{equation}

Let $\varphi(\by,\bi,\bj):= |\Pi_\by(\bi) - \Pi_\by(\bj)|$ and
\begin{equation} \label{def-gr}
g_r(\bi,\bj) := \eta_\infty \{\by\in F_\eps:\ \varphi(\by,\bi,\bj) < r\}.
\end{equation}
Let $P:=G_\varepsilon \times G_\varepsilon $ and
$\mu_2:= \mu_\varepsilon \times \mu_\varepsilon $. Denote $P^N :=
\{(\bi,\bj)\in P:\ |\bi\wedge \bj|\ge N\}$.
  
\begin{prop} \label{prop}
There exists $C=C(\eps)>0$, such that for all $r>0$, 
\begin{equation}\label{19}
A(r):=\iint\limits_{P^N}g_r(\mathbf{i},\mathbf{j})
d\mu_2(\mathbf{i},\mathbf{j})\leq Cr.
\end{equation}
\end{prop}

We will prove the proposition in the next section.
Before that, using the proposition we prove 
Theorem~\ref{main}(a).

\begin{proof}[Conclusion of the Proof of Theorem~\ref{main}(a)] 
In order to prove
(\ref{eq-main2}), it is enough to verify that
$$
\mathcal{I}:=\sum_{|\om|=N}\ \ \int\limits_{\mathbf{F_\eps }}%
\int\limits_{\mathbb{R}}%
\underline{D}(\nu ^\varepsilon  _{\by,\om},x)d\nu^\varepsilon
_{\by,\om}(x)d\eta _\infty (\mathbf{y})<\infty ,
$$
where
$$
\underline{D}(\nu,x):=\liminf_{r\to 0}%
\frac{\nu
\left(\left[x-r,x+r\right]\right)}{2r}\,,
$$
is the lower derivative of a measure $\nu$, see \cite[2.12]{mattila}.
Observe that
$$
\sum_{|\om|=N}\int\limits_{\mathbb{R}}\nu^\varepsilon _{\by,\om}\left(\left[
x-r,x+r\right]\right)d\nu^\eps _{\by,\om}(x)=%
\iint\limits_{P^N}\mathds{1}_{\left\{%
(\mathbf{i},\mathbf{j}):\left|\Pi_\mathbf{y}(\mathbf{i})-%
\Pi_\mathbf{y}(\mathbf{j}) \right|\leq r\right\}}d\mu
_2(\mathbf{i},\mathbf{j})
$$
by the definition of $\nu^\eps _{\by,\om}$.
Using this with Fatou's Lemma, and exchanging the order of integration,
we obtain that
\begin{equation}\label{21}
\mathcal{I }\leq\liminf_{r\to 0}(2r)^{-1}A(r),
\end{equation}
where $A(r)$ was defined in (\ref{19}). Thus
$\mathcal{I}<\infty $ follows immediately from Proposition~\ref{prop}.
\end{proof}


\section{Transversality lemma and the proof of Proposition~\ref{prop}}
                                                                                
We begin with a technical lemma, which is a key for the proof of both parts
of Theorem~\ref{main}. We are assuming all the conditions of
Theorem~\ref{main}, in particular, that either (\ref{cond1}) or
(\ref{cond2}) holds. 

By the definition of $\varphi$ and $\Pi$ we have
\begin{equation}\label{eq-phi}
\varphi (\mathbf{y},\mathbf{i},\mathbf{j}) =
|d_{i_1}-d_{j_1} + y_1 \Phi(\by,\bi,\bj)|,
\end{equation}
where
\begin{equation} \label{eq-Phi}
\Phi(\by,\bi,\bj) =  \lambda_{i_1}d_{i_2}- \lambda_{j_1}d_{j_2}
+\sum\limits_{\ell=2}^{\infty }y_{2\cdots \ell}
(\lambda_{i_{1}\dots i_{\ell}}d_{i_{\ell+1}} -
\lambda_{j_{1}\dots j_{\ell}}d_{j_{\ell+1}}).
\end{equation}
Note that $\Phi(\by,\bi,\bj)$ does not depend on $y_1$.
If $|\mathbf{i}\wedge\mathbf{j}|= k\ge 1$ then
\begin{equation} \label{eq-Phi1}
\varphi (\mathbf{y},\mathbf{i},\mathbf{j})= \lambda _{i_1\dots
i_k}\cdot y_{1\dots k}\cdot \varphi \left(\sig^k\by,\sig^k\bi,\sig^k\bj
\right).
\end{equation}
                                                                                
\begin{lemma} \label{lem1}
Let $\delta>0,\ \rho\in (0,1)$, and $N\in \Nat$. Consider
$$
F = \{\by\in \R^\Nat:\ y_{1\ldots n} \ge (1+\delta)^{-n}, \ \forall\,
n\ge N\},
$$
$$
G = \{\bi \in \Sigma:\ \lam_{i_1\ldots i_n} \ge \rho^n,\ \forall\,
n\ge N\}.
$$
There exists $C_2>0$ such that for all $k\ge N$,
for all $\bi,\bj \in G$, with $|\bi\wedge \bj|=k$,
\begin{equation}\label{eq-r1}
\eta_\infty\{\by \in F: \ \varphi(\by,\bi,\bj) < r\}
 \le C_2 (1+\delta)^k \rho^{-k}  r\ \ \ \mbox{for all}\ r>0.
\end{equation}
\end{lemma}

\begin{proof}
First suppose that the condition (\ref{cond1}) holds. Then
$b:=\min_{\ell\ne s} |d_\ell-d_s|>0$.
Since $\bi \in G$ and $k\ge N$, we have for $\by\in F$ by (\ref{eq-phi}),
(\ref{eq-Phi}) and (\ref{eq-Phi1}):
\begin{equation}\label{eq-Phi3}
\varphi(\by,\bi,\bj) < r
\ \Rightarrow\ |d_{i_{k+1}} - d_{j_{k+1}} + y_{k+1} \Phi|
< (1+\delta)^k \rho^{-k} r,
\end{equation}
where $ \Phi = \Phi (\sig^k\by,\sig^k\bi,\sig^k\bj) $
does not depend on $y_{k+1}$.
Denote $\Delta_{k+1} := d_{i_{k+1}} - d_{j_{k+1}}$; we have
$|\Delta_{k+1}| \ge b$ since $i_{k+1}\ne j_{k+1}$. We can assume that
$\Delta_{k+1}<0$; otherwise, we just switch $\bi$ and $\bj$.
Since the left-hand side of (\ref{eq-r1})
is always bounded above by one, (\ref{eq-r1})
holds for $r\ge (1+\delta)^{-k} \rho^k b/2$ with the constant $C_2=2/b$.
If $r < (1+\delta)^{-k} \rho^k b/2$ then
\begin{equation} \label{eq-Phinew}
\varphi(\by,\bi,\bj) < r
\ \Rightarrow\ |\Delta_{k+1} + y_{k+1} \Phi| < (1+\delta)^k \rho^{-k} r < b/2,
\end{equation}
and this implies $\Phi>0$,
in view of $y_{k+1}$ being positive and the fact that
$\Delta_{k+1}\le -b$. Moreover, the right-hand side of (\ref{eq-Phinew}) implies
$$
y_{k+1}\in B,\ \ \mbox{where}\ \ B:=\left[ 
\frac{-\Delta_{k+1} - (1+\delta)^k\rho^{-k}r}{\Phi}\,,\,
\frac{-\Delta_{k+1} + (1+\delta)^k\rho^{-k}r}{\Phi} \right].
$$
Note that $B$ depends on $y_{k+2},y_{k+3},\ldots$ but not on $y_{k+1}$.
We have $(1+\delta)^k\rho^{-k}r < b/2 \le -\Delta_{k+1}/2$, so
$$
B \subset [-\Delta_{k+1}/(2\Phi), \infty)\subset [b/(2\Phi),\infty).
$$
By (\ref{eq-density}), we obtain that for any $y_{k+2},y_{k+3},\ldots$,
$$
\eta\{y_{k+1}\in B\} \le C_1 (2\Phi/b)\Leb(B) = C_1(4/b) (1+\delta)^k
\rho^{-k} r.
$$
This implies the desired inequality (\ref{eq-r1}) by Fubini Theorem,
since $y_{k+1}$ is independent of $y_{k+2},y_{k+3},\ldots$

\smallskip
                                                                                
Now suppose that the condition (\ref{cond2}) holds. Then
by (\ref{eq-phi}) and (\ref{eq-Phi}),
$$
\varphi(\by,\bi,\bj) = \lam_{i_1\ldots i_k}\cdot y_{1\ldots k,(k+1)}\cdot
|\lam_{i_{k+1}} - \lam_{j_{k+1}} + y_{k+2}\Psi|,
$$
where
$$
\Psi = \lam_{i_{k+1} i_{k+2}}-\lam_{j_{k+1} j_{k+2}}+
\sum_{\ell=3}^\infty y_{(k+3)\ldots (k+\ell)}
(\lam_{i_{k+1}\ldots i_{k+\ell}} - \lam_{j_{k+1}\ldots j_{k+\ell}})
$$
does not depend on $y_{k+2}$. Since $\bi \in G$ and $k\ge N$, we have for
$\by\in F$:
$$
\varphi(\by,\bi,\bj) < r
\ \Rightarrow\ |\lam_{i_{k+1}} - \lam_{j_{k+1}} + y_{k+2} \Psi|
< (1+\delta)^{k+1} \rho^{-k} r.
$$
Here $|\lam_{i_{k+1}} - \lam_{j_{k+1}}|\ge b':=
\min_{\ell\ne s} |\lam_\ell-\lam_s|>0 $ by (\ref{cond2}), and we 
argue similarly to the first case to obtain
(\ref{eq-r1}), with
$C_2 = (1+\delta)\max\left\{\frac{2}{b'}, \frac{4C_1}{b'}\right\}$.
\end{proof}

\begin{proof}[Proof of Proposition \ref{prop}]
Let
$$
P_\omega :=\left\{(\mathbf{i},\mathbf{j})\in P:%
\mathbf{i}\wedge\mathbf{j}=\omega  \right\},
$$
where $\omega=(\omega _1,\dots,\omega _k )\in
\left\{1,\dots,m\right\}^k$ for some $k$.
Denote
$$
A_\omega(r) :=\iint\limits_{P_\omega }g_r(\mathbf{i},\mathbf{j})
d\mu_2(\mathbf{i},\mathbf{j}),
$$
so that
\begin{equation} \label{Ar}
A(r) = \sum_{k=N}^\infty \sum_{|\omega|=k} A_\omega(r).
\end{equation}
We can apply Lemma~\ref{lem1} with 
$N=N(\eps)$. Then $F_\eps \subset F$ and $G_\eps \subset G$ by
(\ref{trunc1}) and (\ref{trunc4}), so for $\bi,\bj\in G_\eps$, with
$|\bi\wedge \bj| = k\ge N$, we have
$$
g_r(\bi,\bj) = \eta_\infty\{\by \in F_\eps:\ \varphi(\by,\bi,\bj) < r\} \le
C_2 (1+\delta)^k \rho^{-k} r.
$$
Thus for $|\om|=k\ge N$,
\begin{equation} \label{last1}
A_\om(r) \le C_2 (1+\delta)^k \rho^{-k} r
\cdot (\mu\times\mu)\{(\bi,\bj):\ \bi \wedge \bj = \om\}.
\end{equation}
On the other hand,
\begin{equation} \label{eq-now}
(\mu\times\mu)\{(\bi,\bj):\ \bi \wedge \bj = \om\} \le \mu([\om])^2 =
\mu[\bi,k]\cdot \mu([\om]) \le \theta^k \mu([\om]),
\end{equation}
in view of (\ref{trunc4}).
Combining this with (\ref{last1}) and (\ref{Ar}) we obtain
$$
A(r) \le C_2 \sum_{k\ge N} \sum_{|\om|=k} (1+\delta)^k
\theta^k\rho^{-k} \mu([\om]) \cdot r < \const \cdot r,
$$
where we used that $\sum_{|\om|=k}\mu([\om])=1$ and
(\ref{delta}). The proof is complete.
\end{proof}


\section{Proof of Theorem~\ref{main}(b)}

Fix $\eps\in (0,1)$; our goal is to prove (\ref{eq-main3}) assuming that
$-h(\mu) \ge  \chi(\mu)$.

\smallskip

{\sc Estimate from below.} Fix an arbitrary $\alpha < h(\mu)/|\chi(\mu)|$;
it is enough to prove that 
\begin{equation} \label{eq-below}
\Dh (\nu_\by^\eps) \ge \alpha \ \ \mbox{for}\
\eta_\infty\ \mbox{a.e.\ } \by \in F_\eps.
\end{equation}
We can find $\theta,\rho,\delta>0$ such that
\begin{equation} \label{eq-condit}
\alpha < \frac{\log \theta}{\log((1+\delta)^{-1}\rho)},\ \ \
|\chi(\mu)| < - \log\rho,\ \ \mbox{and}\ \ h(\mu) > - \log\theta.
\end{equation}
Similarly to the proof of Theorem~\ref{main}(a), 
we can find $N = N(\eps,\delta,\theta,\rho)$ such that 
$$
(1+\delta)^{-n} \le y_1\cdots y_n \le (1+\delta)^n\ \ \ \mbox{for all}\ \
n\ge N,\ \by \in F_\eps
$$
and
$$
\mu[\bi,n]  < \theta^n\ \ \ \mbox{and}\ \ \ \rho^n  < \lam_{i_1\ldots i_n},
\ \ \ \mbox{for all}\ \
n\ge N,\ \bi \in G_\eps.
$$
We will use the decomposition (\ref{decomp}) again.

By Frostman's Theorem, see \cite[Theorem 4.13]{Falc1}, for any Borel
measure $\nu$ on the line,
\begin{equation} \label{eq-frost}
\Dh(\nu) \ge \sup\lt\{\alpha>0:\ \iint\limits_{\R^2}
\frac{d\nu(\xi)\,d\nu(\zeta)}{|\xi-\zeta|^\alpha}
<\infty\rt\}\,.
\end{equation}
Thus the desired estimate (\ref{eq-below}) will follow by Fubini's Theorem,
if we show that
\begin{equation} \label{eq-Sk}
\Sk:= \sum_{|\om|=N} \int\limits_{F_\eps} \iint\limits_{\R^2}
|\xi-\zeta|^{-\alpha}\,d\nu^\eps_{\by,\om}(\xi)\,d\nu^\eps_{\by,\om}(\zeta)
\, d\eta_\infty(\by) < \infty.
\end{equation}
After changing the variables and reversing the order of integration we
obtain
\begin{equation} \label{eq-Sk1}
\Sk = \sum_{k=N}^\infty \sum_{|\om|=k} \iint\limits_{P_\omega }
\int\limits_{F_\eps} 
\varphi(\by,\bi,\bj)^{-\alpha}\,d\eta_\infty(\by)\,d\mu_2(\bi,\bj),
\end{equation}
where again
$$
P_\omega = \{(\bi,\bj) \in G_\eps \times G_\eps:\ \bi\wedge \bj = \om\}.
$$
Suppose that $|\bi\wedge \bj| = k$.
The inner integral in (\ref{eq-Sk1}) is equal to
$$
\alpha \int_0^\infty \eta_\infty \{\by\in F_\eps:\ \varphi(\by,\bi,\bj) \le r\}
\,r^{-\alpha-1}\,dr = \int_0^{(1+\delta)^{-k} \rho^k}  + 
\int_{(1+\delta)^{-k} \rho^k}^\infty .
$$
The first integral in the right-hand side
is estimated by Lemma~\ref{lem1}, and the
second integral is estimated by the trivial estimate
$\eta\{\cdot\} \le 1$ yielding the inequality
$$
\int_{F_\eps} \varphi(\by,\bi,\bj)^{-\alpha}\,d\eta_\infty(\by) \le
\const \cdot [(1+\delta)\rho^{-1}]^{\alpha k}.
$$
Substituting this into (\ref{eq-Sk1}) we obtain
$$
\Sk \le \const\cdot \sum_{k=N}^\infty \sum_{|\om|=k}
[(1+\delta)\rho^{-1}]^{\alpha k} \mu_2\{(\bi,\bj):\ \bi \wedge \bj = \om\}.
$$
Now we can apply (\ref{eq-now}) to get
$$
\Sk \le \const\cdot \sum_{k=N} [(1+\delta)\rho^{-1}]^{\alpha k} \theta^k <
\infty,
$$
by
(\ref{eq-condit}).

\medskip

{\sc Estimate from above.} Dimension estimates from above are 
fairly standard.
This is also the case here, although there are
technical complications because of the generality of our set-up.
Note that we obtain the upper
bound for all, rather than almost all, $\by \in F_\eps$, and the distribution of
$y_i$'s is irrelevant here. (Recall that $F_\eps$ was defined at the
beginning of Section 3.) A similar upper bound for 
(possibly nonlinear, but non-random)
contracting on average IFS was obtained in \cite{NSB,FST}.

Fix an arbitrary $\alpha > h(\mu)/|\chi(\mu)|$;
it is enough to prove that
\begin{equation} \label{eq-above}
\forall\,\eps>0,\
\Dh (\nu_\by) \le \alpha \ \ \mbox{for all}\ \by \in F_\eps.
\end{equation}
We fix $\eps>0$ and $\by \in F_\eps$ for the rest of this proof.

Now let $\gam>0$ and consider $G_\gam$, with
$\mu(G_\gam) > 1-\gam$, such that the convergence in (\ref{shannon}) and
(\ref{ergth}) is uniform on $G_\gam$.
Further, let $C_3>0$ be such that
\begin{equation} \label{eq-Omega}
\mu(\Omega_\gam) \ge 1-\gam,\ \ \ \mbox{where}\ \ 
\Omega_\gam:= \{\bi \in \Sigma:\ |\Pi_\by(\bi)| \le C_3\}.
\end{equation}
Consider
\begin{equation} \label{eq-A}
A_\gam^n:= \left\{\bi \in G_\gam \cap \sig^{-n} \Omega_\gam:\ 
\mu([\bi,n]\cap \sig^{-n} \Omega_\gam) \ge 0.5\cdot \mu([\bi,n])\right\}\,.
\end{equation}
We claim that $\mu(A_\gam^n) \ge 1-4\gam$ for all $n\in \Nat$. Indeed,
$\mu(G_\gam \cap \sig^{-n} \Omega_\gam) \ge 1-2\gam$, and the measure of the
complement of the set of $\bi$ satisfying the inequality in 
(\ref{eq-A}) equals
$$
\mu\left\{\bi \in \Sigma:\ 
\mu([\bi,n]\cap (\sig^{-n} \Omega_\gam)^c) \ge 0.5\cdot \mu([\bi,n]\right\}
\le 2 \mu((\sig^{-n} \Omega_\gam)^c) \le 2\gam,
$$
where we used that $\mu$ is $\sigma$-invariant in the last step.
It follows that
\begin{equation} \label{eq-H}
\mu(H_\gam) \ge 1-4\gam,\ \ \ \mbox{where}\ \ 
H_\gam:= \limsup(A_\gam^n)=\bigcap_{n=1}^\infty \bigcup_{k=n}^\infty A_\gam^n.
\end{equation}
Recall that our goal is to prove $\Dh (\nu_\by) \le \alpha$. 
Billingsley's Theorem (see \cite[p.171]{falcbook3}) states that
$$
\Dh (\nu_\by)=\nu_\by\mbox{-}\ess \sup\left\{\liminf_{r\downarrow 0}
\frac{\log \nu_\by[x-r,x+r]}{\log (2r)}\right\}.
$$
Thus it is enough to verify that
\begin{equation} \label{eq-bil}
\liminf_{r\downarrow 0}
\frac{\log \nu_\by[x-r,x+r]}{\log (2r)} \le \alpha
\end{equation}
for $\nu_\by$ a.e.\ $x$.
Since $\nu_\by = \mu \circ \Pi_\by^{-1}$ and in view of (\ref{eq-H}), 
this will follow if we prove (\ref{eq-bil}) for all $x\in \Pi_\by(H_\gam)$,
for every $\gam>0$. To this end, let us fix $\gam>0$ and $x = \Pi_\by(\bi)$
for some $\bi \in H_\gam$. 
Since $\bi \in H_\gam$, there exists a sequence $n_k \to \infty$ such that
\begin{equation} \label{eq-esti}
\mu([\bi,n_k] \cap \sigma^{-n_k} \Omega_\gam) \ge 0.5\cdot \mu[\bi,n_k],\ \ \ \
\forall\, k\in \Nat.
\end{equation}

Since $\alpha > h(\mu)/|\chi(\mu)|$, we 
can find $\theta,\rho,\delta>0$ such that
\begin{equation} \label{eq-condit2}
\alpha > \frac{\log \theta}{\log((1+\delta)\rho)}\,,\ \ \
|\chi(\mu)| >  - \log\rho,\ \ \mbox{and}\ \ h(\mu) < - \log\theta.
\end{equation}
Similarly to the proof of Theorem~\ref{main}(a),
we can find $N$ such that (\ref{trunc1}) holds and
\begin{equation} \label{trunc44}
\mu[\bi,n]  > \theta^n \ \ \ \mbox{and}\ \ \  \rho^n  > \lam_{i_1\ldots i_n},
\ \ \ \mbox{for all}\ \
n\ge N,\ \bi \in G_\gam.
\end{equation}
Let $r_k = 2C_3 \rho^{n_k} (1+\delta)^{n_k}$. We claim that for all $k$
sufficiently large,
\begin{equation} \label{eq-estim}
\nu_\by[x-r_k,x+r_k] = \mu\{\bj:\ |\Pi_\by(\bi)-\Pi_\by(\bj)| \le r_k\}
\ge 0.5\cdot \mu[\bi, n_k]
\end{equation}
(the equality here is by definition; the claim is the inequality).
Indeed, let $\bj \in [\bi,n_k] \cap \sigma^{-n_k} \Omega_\gam$. Then for
$k$ sufficiently large (so that $n_k \ge N$), we have
$$
|\Pi_\by(\bi)-\Pi_\by(\bj)| = \lam_{i_1\ldots i_{n_k}} y_{1\ldots n}
\cdot |\Pi_\by(\sig^{n_k}\bi)-\Pi_\by(\sig^{n_k}\bj)| \le
2C_3 \rho^{n_k} (1+\delta)^{n_k}=r_k,
$$
using (\ref{trunc44}), (\ref{trunc1})
and the fact that
$\sig^{n_k}\bi, \sig^{n_k}\bj\in \Omega_\gam$,
where $\Omega_\gam$ is defined by (\ref{eq-Omega}).
This, combined with (\ref{eq-esti}), proves (\ref{eq-estim}). Now,
keeping in mind that the numerator and denominator below are negative,
we obtain
\begin{eqnarray*}
\liminf_{k\to\infty}
\frac{\log \nu[x-r_k,x+r_k]}{\log (2r_k)} & \le & 
\liminf_{k\to \infty} 
\frac{\log \mu[\bi, n_k]-\log 2}{\log(2C_3) + n_k\log((1+\delta)\rho)} \\
& \le & \lim_{k\to \infty}
\frac{n_k\log\theta}{n_k \log((1+\delta)\rho)} < \alpha,
\end{eqnarray*}
where we used  (\ref{trunc44}) and
(\ref{eq-condit2}). The proof is complete. \qed
 

\section{Proofs of other results}

\subsection{Method of Kahane-Salem}
Here we prove Theorem~\ref{th-homog} using a variant of the approach from
\cite{KS}.
Recall that for a finite measure $\nu$ on $\R$ its Fourier transform is
defined by $\nuhat(\xi) = \int_\R e^{it\xi}\,d\nu(t)$.

\begin{definition} \label{def-sobolev}
{\em 
For a finite measure $\nu$ on $\R$, its {\em Sobolev dimension} is defined as
\begin{equation} \label{eq-sobol}
\dim_s(\nu) = \sup\left\{\alpha \in \R:\ \Ek_\alpha(\nu)=\int_{\R}
|\nuhat(\xi)|^2(1+|\xi|)^{\alpha-1}\,d\xi < \infty\right\}\,.
\end{equation}}
\end{definition}
                                                                                
\begin{remark} \label{rem6}
{\em If $\dim_s(\nu) < 1$, then $\dim_s(\nu)$ is also known as the correlation
dimension of the measure $\nu$. If $\Ek_\alpha(\nu)< \infty$ for
$\alpha>1$, then $\nu \ll \Leb$,
and its density is said to have the fractional derivative of order 
$(\alpha-1)/2$
in $L^2(\R)$. If $\Ek_1(\nu)<\infty$, then $\nu$ has a
density in $L^2(\R)$ (this is just Plancherel's Theorem), 
and if $\dim_s(\nu)>2$, then $\nu$ has a
continuous density, see e.g.\ \cite[Th.\ 1.2.4]{AH}.}
\end{remark}

\begin{theorem} \label{th-Sobo}
Let $\nu_\by$ be the conditional distribution of the sum (\ref{eq-proj})
given $\by=(y_1,y_2,\ldots)$, defined by (\ref{eq-nuby}). We assume that
(\ref{eq-dens2}) and (\ref{cond1}) are
satisfied, $\mu = (p_1,\ldots,p_m)^\Nat$ and $\lam_i=\lam$ for all $i\le m$.
Suppose that $\chi(\eta) = \log\lam+\E[\log Y]<0$ and denote
$\beta = \sum_{j=1}^m p_j^2$.
Then
$$
\dim_s(\nu_\by) \ge |\log\beta|/|\chi(\eta)|\ \ \ 
\mbox{for}\ \ \eta_\infty\ a.e.\ \by.
$$
In particular, if $0>\chi(\eta)>\log\beta$, then $\nu_\by \ll \Leb$ with
a density in $L^2(\R)$ for $\eta_\infty$ a.e.\ $\by$. 
If $0> \chi(\eta) > \half\log\beta$, then $\nu_\by \ll \Leb$ with
a continuous density for $\eta_\infty$ a.e.\ $\by$.
\end{theorem}

In view of Remark~\ref{rem6},
Theorem~\ref{th-homog} is contained in Theorem~\ref{th-Sobo}.

\medskip

\noindent {\em Proof of Theorem~\ref{th-Sobo}.} 
Since $\lam_i=\lam$ for all $i\le m$, we can assume without loss of
generality that $\lam=1$ (just replace $Y$ with $\lam Y$). Then $\chi(\eta) = 
\E[\log Y]$.
Our goal is to prove that for every $\alpha < 
\frac{|\log\beta|}{|\chi(\eta)|}$,
\begin{equation}\label{eq1}
\int_\R |\nuhat_\by(\xi)|^2(1+|\xi|)^{\alpha-1}\,d\xi < \infty,
\end{equation}
for $\eta_\infty$ a.e.\ $\by$.
Fix $\alpha < \frac{|\log\beta|}{|\chi(\eta)|}$ for the rest of the proof.
By the the Law of Large Numbers and Egorov's Theorem, 
for any $\eps>0$ we can find $F_\eps\subset \R^\Nat$ such that
$\eta_\infty(F_\eps) > 1-\eps$ and
$(y_{1\ldots n})^{1/n} \to e^{\chi(\eta)}$ uniformly on $F_\eps$.
It suffices to verify (\ref{eq1}) for $\eta_\infty$ a.e.\ $\by \in F_\eps>0$,
for an arbitrary $\eps>0$. Fix $\eps > 0$ for the rest of the proof. The
result will follow by Fubini's Theorem if we can show that
$$
\int_{F_\eps} \int_\R |\nuhat_\by(\xi)|^2(1+|\xi|)^{\alpha-1}\,d\xi\,
d\eta_\infty (\by) < \infty.
$$

Recall that $\nu_\by$ is the conditional distribution of $X$ in (\ref{absorb}) 
given $\by$, with $\lam=1$, which can be viewed as a sum of independent
discrete random variables. Thus, 
$\nu_\by$ is the infinite convolution product
$$
\nu_{\by} = \mbox{\huge $\ast$}\prod_{n=1}^\infty \left( \sum_{j=1}^m p_j
\delta_{d_j y_{1\ldots (n-1)}} \right)\,,
$$
where $\delta$ is the Dirac's delta. Its Fourier transform is
\begin{equation} \label{Fou1}
\nuhat_\by(\xi) 
=\prod_{n=1}^\infty \sum_{j=1}^m p_j e^{id_j y_{1\ldots (n-1)} \xi}
=: \prod_{n=1}^\infty \psi_n(\by,\xi).
\end{equation}
Now the argument essentially follows the proof of \cite[Th\'eor\`eme II]{KS}.
We have
\begin{eqnarray*}
|\psi_n(\by,\xi)|^2 & = & \sum_{j,k=1}^m p_j p_k 
e^{i(d_j-d_k) y_{1\ldots (n-1)} \xi} \\ & = &
\sum_{j=1}^m p_j^2 + \sum_{j\ne k} p_j p_k
e^{i(d_j-d_k) y_{1\ldots (n-1)} \xi}.
\end{eqnarray*}
Clearly,
$$
|\nuhat_\by(\xi)|^2 \le  \prod_{n=1}^{\ell+1} |\psi_n(\by,\xi)|^2 =: 
f_{\xi,\ell+1}.
$$
Denote by $F_\eps^n$ the projection of $F_\eps$ onto $\R^n$ (the first $n$
coordinates). Then, since $f_{\xi,\ell+1}$ depends only on $y_1,\ldots,y_\ell$,
we obtain
\begin{eqnarray} 
\int\limits_{F_\eps} f_{\xi,\ell+1}\,d\eta_{\infty}(\by) & = & 
\int\limits_{F_\eps^\ell} f_{\xi,\ell+1}\,d\eta(y_1)\ldots d\eta(y_\ell)
\nonumber \\
& \le & \int\limits_{F_\eps^{\ell-1}} f_{\xi,\ell} \,d\eta(y_1)\ldots 
d\eta(y_{\ell-1})\,\int\limits_\R |\psi_{\ell+1}(\by,\xi)|^2\,d\eta(y_\ell).
\label{eq2}
\end{eqnarray}
Recall that $\beta = \sum_{j=1}^m p_j^2$, so
\begin{eqnarray}
\int_\R |\psi_{\ell+1}(\by,\xi)|^2\,d\eta(y_\ell) & = &
\beta + \sum_{j\ne k} p_j p_k \int_\R e^{i(d_j-d_k) y_{1\ldots (n-1)} \xi}
\,d\eta(y_\ell) \nonumber \\
& = & \beta + \sum_{j\ne k} p_j p_k 
\widehat{\eta}((d_j-d_k) y_{1\ldots (n-1)} \xi).
\label{eq3}
\end{eqnarray}
Integration by parts (see e.g. \cite[p.\ 25]{Katz}) shows
that the Fourier transform of
 a compactly supported function of bounded variation is bounded above by
$c|t|^{-1}$. Since $|d_j-d_k| \ge b> 0$ for $j\ne k$, we obtain
that for some $C>0$,
\begin{equation} \label{eq4}
\int_\R |\psi_{\ell+1}(\by,\xi)|^2\,d\eta(y_\ell) \le \beta
\left( 1 + \frac{C}{y_{1\ldots (\ell-1)}|\xi|}\right)
\end{equation}
Recall that $ \alpha < |\log \beta|/|\chi(\eta)|$; choose $\rho < 
e^{\chi(\eta)}$ such 
that $\alpha < |\log \beta|/|\chi(\eta)|$. Since $y_{1\ldots n}^{1/n}$
converges to $e^{\chi(\eta)}$ 
uniformly on $F_\eps$, we can find $N\in \Nat$ such that
$$
y_{1\ldots n} \ge (n+1)^2\rho^{n+1}\ \ \ \forall\, n\ge N,\ \forall\,
\by \in F_\eps.
$$
It follows from (\ref{eq2}) and (\ref{eq4}) that for $\ell \ge N+1$,
\begin{eqnarray}
\int_{F_\eps} f_{\xi,\ell+1}\,d\eta_{\infty}(\by) & \le &
\int_{F_\eps} f_{\xi,\ell}\,d\eta_{\infty}(\by) \cdot \beta
\left( 1 + \frac{C}{\ell^2\rho^\ell |\xi|}\right) \nonumber \\
& \le & \int_{F_\eps} f_{\xi,\ell}\,d\eta_{\infty}(\by) \cdot \beta
(1+C\ell^{-2}), \label{eq5}
\end{eqnarray}
provided that $|\xi| \ge \rho^{-\ell}$. Clearly this condition holds
for $\ell'<\ell$ if it holds for $\ell$, so we can iterate (\ref{eq5})
to obtain, assuming  $|\xi| \ge \rho^{-\ell}$:
$$
\int_{F_\eps} f_{\xi,\ell+1}\,d\eta_{\infty}(\by) \le
\int_{F_\eps} f_{\xi,N}\,d\eta_{\infty}(\by)\cdot \beta^{\ell+1-N}
\prod_{k=N}^\ell(1+C k^{-2}) \le C' \beta^\ell,
$$
where $C'>0$ depends on $N$ but not on $\xi$. We conclude that
$$ 
\int_{F_\eps} |\nuhat_\by(\xi)|^2\,d\eta_\infty(\by) \le C''
|\xi|^{\log\beta/|\log\rho|}.
$$
But $\log\beta/|\log\rho| < -\alpha$, hence
$$
\int_\R \left( \int_{F_\eps} |\nuhat_\by(\xi)|^2\,d\eta_\infty(\by) \right)
(1+|\xi|)^{\alpha-1}\,d\xi < \infty,
$$
and the proof is complete. \qed

\subsection{Proof of Proposition~\ref{prop-sing}}

Consider the probability space $\Sig \times \R^\Nat$ with the measure $\Prob:=
\mu\times \eta_\infty$.
Under the assumptions of the proposition, $Z_n: =\log p_{i_n} - \log \lam_{i_n}
-\log y_n$ are i.i.d.\ non-constant random variables with mean zero.
Let
$$
B_n := \left\{(\bi,\by):\ \sum_{j=1}^n Z_j > \sqrt{n} \right\}
= \left\{ (\bi,\by):\ 
\frac{p_{i_1\ldots i_n}}{\lam_{i_1\ldots i_n} y_{1\ldots n}} > 
e^{\sqrt{n}}\right\}.
$$
By the Law of Iterated Logarithm, we have 
\begin{equation} \label{eq-limsup}
\Prob(\limsup B_n) =1.
\end{equation}
For  $(\bi,\by)\in \limsup B_n$ and $k\in \Nat$ let
$$
\tau_k=\tau_k(\bi,\by) = \min\{n\ge k:\ (\bi,\by)\in B_n\},
$$
which is well-defined and finite. Let
$$
A_n = \{(\bi,\by):\ |\Pi_{\sig^n\by}(\sig^n\bi)|\le n\};
$$
note that $A_n$ is independent of $i_1,\ldots,i_n, y_1,\ldots, y_n$.
Since the probability $\Prob$ is $\sig$-invariant and
$\Pi_\by(\bi)$ is finite a.s., we have $\Prob(A_n)=\Prob\{(\bi,\by):\
\Pi_\by(\bi) \le n\} \uparrow 1$. Now,
\begin{eqnarray*}
\Prob\left[ \bigcup_{n\ge k} (A_n \cap B_n)\right] & \ge & 
\Prob\left( \bigcup_{n\ge k} [A_n \cap\{\tau_k=n\}]\right) \\
& = & \sum_{n=k}^\infty \Prob\{\tau_k=n\}\cdot \Prob(A_n) \ge \Prob(A_k) \to 1,
\end{eqnarray*}
as $k\to \infty$. In the last displayed line we used that 
$A_n$ and $\{\tau_k=n\}$ are independent events and that
$\sum_{n=k}^\infty \Prob\{\tau_k=n\}=1$ by (\ref{eq-limsup}). It follows that
$$
\Prob(\limsup(A_n\cap B_n)) =1.
$$
By Fubini, there exists $\Omega\subset \R^\Nat$ 
such that $\eta_\infty(\Omega)=1$ and 
$$
\Sig_\by:= \{\bi\in \Sig:\ (\bi,\by)\in \limsup(A_n\cap B_n)\}
$$
has $\mu(\Sig_\by) =1$
for every $\by\in \Omega$. We claim that $\Leb(\Pi_\by(\Sig_\by))=0$ for
every $\by\in \Omega$, which will imply that $\nu_\by =\mu\circ \Pi_\by^{-1}
\perp\Leb$.
Fix $\by\in \Omega$ and $k \in \Nat$. 
Observe that
$\Sig_\by \subset \bigcup_{n\ge k} \{\bi: (\bi,\by)\in A_n\cap B_n\}$.
For any $(\bi,\by)\in A_n\cap B_n$ and $(\bj,\by) \in [\bi,n]\cap A_n$ we have
\begin{eqnarray*}
|\Pi_\by(\bi)-\Pi_\by(\bj)| & = & \lam_{i_1\ldots i_n}y_{1\ldots n} 
|\Pi_{\sig^n\by}(\sig^n\bi) - \Pi_{\sig^n\by}(\sig^n\bj)| \\
& \le & 2n \lam_{i_1\ldots i_n}y_{1\ldots n} \\
& \le & 2n e^{-\sqrt{n}} p_{i_1\ldots i_n}. 
\end{eqnarray*}
Here we used first that $(\bi,\by),(\bj,\by)\in A_n$, 
and then that $(\bi,\by)\in B_n$.
Summing over all cylinders of length $n$ (using that $\sum_{i_1\ldots i_n} 
p_{i_1\ldots i_n} =1$), and then summing over $n$ we obtain
$$
\Leb(\Pi_\by(\Sig_\by))  \le \sum_{n\ge k} 2ne^{-\sqrt{n}} \to 0,\mbox{\ \ as\ }
k\to \infty.
$$
\qed

\subsection{Proof of Corollaries~\ref{cor-ranset}-\ref{cor-zero}}
By the Variational Principle (see e.g. \cite{Walters}), $h_{\rm top}(\Gam)
= \sup_\mu h(\mu)$, where the supremum is over ergodic $\sigma$-invariant
measures supported on $\Gam$. Thus, Theorem~\ref{main}(a) implies
Corollary~\ref{cor-ranset}(a), and Theorem~\ref{main}(b) implies 
the lower estimate for $\Dh(S_\Gam(\by))$ in Corollary~\ref{cor-ranset}(b).

In Corollary~\ref{cor-ranset1}, we have $\Gam=\Sig$, for which the
measure of maximal entropy is $(\frac{1}{m},\ldots,\frac{1}{m})^\Nat$.
Part (a) of Corollary~\ref{cor-ranset1} then follows from
Theorem~\ref{th-homog}(b) by the Variational Principle, and 
Corollary~\ref{cor-ranset1}(b) is a special case of 
Corollary~\ref{cor-ranset}(b).

It remains to verify the upper estimate for 
$\Dh(S_\by)$ in Corollary~\ref{cor-ranset}(b).
By the Law of Large Numbers and Egorov's Theorem, we can find $F_\eps\subset
\R^\Nat$ such that $\eta_\infty(F_\eps)>1-\eps$ and $y_{1\ldots n}^{1/n}
\to e^{\chi(\eta)}$ uniformly for $\by\in F_\eps$.
Fix an arbitrary $\alpha > h_{\rm top}(\Gam)/|\chi(\eta)|$. 
It suffices to show that
for every $\eps>0$ we have $\udim_{\rm box}(S_\by)<\alpha$, 
for a.e.\ $\by\in F_\eps$. Here $\udim_{\rm box}$ 
denotes the upper box-counting dimension.
Fix $\eps>0$ for the rest of the proof.

Let $\delta>0$ be such that $\chi(\eta)+\delta<0$ and $\alpha > 
h_{\rm top}(\Gam)/(|\chi(\eta)|-\delta)$.
We can find $N\in \Nat$ such that
$$
y_{1\ldots n} \le \exp(n(\chi(\eta)+\delta)),\ \ \forall\, n\ge N,\ \forall\,
\by\in F_\eps.
$$
Then for $\bi,\bj\in \Sig$ such that $|\bi\wedge\bj|\ge n \ge N$, we have
\begin{eqnarray*}
|\Pi_\by(\bi)-\Pi_\by(\bj)| & = & \left|\sum_{\ell=n}^\infty y_{1\ldots \ell}
(d_{i_{\ell+1}}-d_{j_{\ell+1}}) \right| \\
& \le & 2 d_{\max} \sum_{\ell=n}^\infty e^{\ell(\chi(\eta)+\delta)}
= \frac{2d_{\max} e^{n(\chi(\eta)+\delta)}}{1-e^{\chi(\eta)+\delta}},
\end{eqnarray*}
where $d_{\max}:= \max_{i\le m}|d_i|$. It follows that the diameter of the
set $\Pi_\by([\om])$ for $\om$ of length $n\ge N$ and $\by\in F_\eps$ is
bounded above by $\const\cdot e^{n(\chi(\eta)+\delta)}$. Denote by 
$\#\Wk_n(\Gam)$
the number of cylinders $[\om]$ of length $n$ such that
$[\om]\cap \Gam \ne \es$. We obtain a cover of $S_\by$ by
$\#\Wk_n(\Gam)$ 
intervals of length $\const\cdot e^{n(\chi(\eta)+\delta)}$, hence
$$
\udim_{\rm box}(S_\by) \le \limsup_{n\to \infty}
\frac{\log(\#\Wk_n(\Gam))}{n(-\chi(\eta)-\delta)} = 
\frac{h_{\rm top}(\Gam)}{|\chi(\eta)|-\delta} < \alpha.
$$
Here we used the definition of the upper box-counting dimension and
the definition of topological entropy. The proof is complete. \qed

\medskip

\noindent {\em Proof of Corollary~\ref{cor-zero}.}
This essentially follows the proof of Proposition~\ref{prop-sing}.
Let $\mu = (\frac{1}{m},\ldots \frac{1}{m})$. Then
$Z_n = -\log m - \log \lam - \log y_n$ are i.i.d.\ non-constant
random variables with mean zero on $\R^\Nat$.
Define $B_n' = \left\{\by:\ \sum_{j=1}^n Z_j > \sqrt{n} \right\}.$
Then $\eta_\infty(\limsup B'_n) =1$ by the Law of Iterated Logarithm. We can define $\tau_k=\tau_k(\by)$ similarly to 
the proof of Proposition~\ref{prop-sing}. Then let
$A_n' = \{\by:\ |\Pi_{\sig^n\by}(\sig^n\bi)|\le n\ \forall\,
\bi\in \Sig\}$. We have $\eta_\infty(A'_n)\uparrow 1$ and
$\eta_\infty(\limsup(A'_n\cap B'_n)) =1$ repeating the argument in 
the proof of Proposition~\ref{prop-sing}. 
Now we can take $\Sig_\by = \Sig$ and conclude as in 
the proof of Proposition~\ref{prop-sing}, obtaining that
$\Leb(S_\by) = \Leb(\Pi_\by(\Sig)) =0$ for all $\by \in
\limsup(A'_n\cap B'_n)$. \qed

                                                                                
\section{Open questions}

\noindent {\bf Question 1.} Is the condition  
$\log m > 2|\chi(\eta)|$ in 
Corollary~\ref{cor-ranset1}(a) for the random set $S_\by$ to contain an
interval almost surely, sharp? Perhaps, $\log m > |\chi(\eta)|$ is
already sufficient? This would mean that as soon the random set  $S_\by$
has positive Lebesgue measure, it has non-empty interior (almost surely). 
This is interesting, in particular,
for the example considered 
by Arratia, see Proposition~\ref{prop-ex2}. 

\medskip

\noindent {\bf Question 2.} 
The results on continuous density and intervals in
random sets (see Theorem~\ref{th-homog} and Corollary~\ref{cor-ranset1}(a))
are obtained only in the case when $\mu$ is a product measure.
Extend this to the general  case of ergodic $\mu$.

\medskip 

\noindent {\bf Question 3.} In Theorem~\ref{th-Sobo} we prove
a lower bound for the a.s.\ value of the Sobolev dimension
$\dim_s(\nu_\by)$. Is this actually an equality? If $|\log\beta| < 
|\chi(\eta)|$, then
the matching upper bound can be obtained from the fact that the
Sobolev dimension equals the correlation dimension when it is less than one,
but what about the case $|\log\beta| \ge |\chi(\eta)|$?

\medskip

\noindent{\bf Acknowledgment.} We are grateful to Richard Arratia and Jim
Pitman for useful discussions.


\end{document}